\documentclass{llncs}

\usepackage{latexsym}
\usepackage{amsfonts}
\usepackage{graphicx}
\usepackage{amssymb}
\usepackage{amssymb}

\title{An authentication scheme based on the twisted conjugacy problem}

\author{Vladimir Shpilrain \and Alexander Ushakov}
\institute{Department of Mathematics, The City  College  of New York, 
NY 10031, USA \email{shpil@groups.sci.ccny.cuny.edu} 
\thanks{Research of the first author was partially supported by
the NSF grant DMS-0405105.}
 \and Department of Mathematics, Stevens Institute
of Technology, Hoboken, NJ 07030, USA \email{aushakov@stevens.edu} }

\newcommand{\mod}{~mod~}

\begin{document}

\maketitle
\begin{abstract}
The conjugacy search problem in a group $G$ is the problem of
recovering an $x \in G$ from given $g \in G$ and $h=x^{-1}gx$. The
alleged computational hardness of this problem  in some groups was
used in several recently suggested public key exchange protocols,
including the one due to Anshel, Anshel, and Goldfeld, and the one
due to Ko, Lee et al. Sibert,  Dehornoy,  and  Girault used this
problem in their authentication scheme, which was inspired by the
Fiat-Shamir scheme  involving repeating several times a three-pass
challenge-response step.

In this paper, we offer an authentication scheme whose security is
based on the apparent hardness of the {\it twisted conjugacy search
problem} which is: given a pair of endomorphisms (i.e.,
homomorphisms into itself) $\varphi, \psi$ of a group $G$ and a pair
of elements $w, t \in G$, find  an element $s \in G$ such that $t =
\psi(s^{-1}) w \varphi(s)$ provided   at least one such $s$ exists.
This problem appears to be very non-trivial even for free groups. We
offer here another platform, namely, the  {\it semigroup} of all $2
\times 2$ matrices over truncated one-variable polynomials over
${\mathbf F}_2$, the field of two elements, with transposition used
instead of inversion in the equality above.
\end{abstract}

\section{Introduction}

One of the most obvious ramifications of the discrete logarithm
problem in the noncommutative situation is the {\it conjugacy search
problem}:

\begin{quote}
Given a  group $G$ and two conjugate elements $g, h \in G$, find a
particular element $x \in G$ such that $x^{-1}gx=h$.
\end{quote}

 This problem always has a recursive solution because one can recursively enumerate all
conjugates of a given element, but this kind of solution can be
extremely inefficient. Specific groups may or may not admit more
efficient solutions, so the choice of the platform group is of
paramount importance for security of a cryptographic primitive based
on the conjugacy search problem. A great deal of  research was (and
still is) concerned with the complexity of this problem in braid
groups because there were several proposals, including the one by
Anshel, Anshel, and Goldfeld \cite{AAG}, and the one by Ko, Lee at
al. \cite{KLCHKP} on using the alleged computational hardness of
this problem  in braid groups to build a key exchange protocol.
Also, Sibert,  Dehornoy,  and  Girault  \cite{SDG} used this problem
in their authentication scheme, which was inspired by the
Fiat-Shamir scheme involving repeating several times a three-pass
challenge-response step. At the time of this writing, no
deterministic polynomial-time algorithm for solving the conjugacy
search problem in braid groups has been reported yet; see
\cite{BGM1} and \cite{BGM2} for recent progress in this direction.
However, several heuristic algorithms, in particular so-called
``length based attacks", were shown to have very high success rates,
see e.g. \cite{GKTTV}, \cite{GKTTV2}, \cite{HS}, \cite{MU3},
\cite{Ruinskiy}. This shows that one has to be really careful when
choosing the platform (semi)group to try to avoid length based or
similar attacks. One way to achieve this goal is, informally
speaking, to have ``a lot of commutativity" inside otherwise
non-commutative (semi)group; see \cite{Ruinskiy} for a more detailed
discussion.

In this paper, we propose an authentication scheme whose security is
based on the apparent hardness of the {\it (double) twisted
conjugacy search problem} which is:

\begin{quote}
given a pair of endomorphisms (i.e., homomorphisms into itself)
$\varphi, \psi$ of a group $G$ and a pair of elements $w, t \in G$,
find  an element $s \in G$ such that $t = \psi(s^{-1}) w \varphi(s)$
provided   at least one such $s$ exists.
\end{quote}

This problem, to the best of our knowledge, has not been considered
in group theory before, and neither was its decision version: given
$\varphi, \psi \in End(G)$, $w, t \in G$, find out whether or not
there is an element $s \in G$ such that $t = \psi(s^{-1}) w
\varphi(s)$. However, the  following special case of this problem
(called the {\it twisted conjugacy problem}) has recently attracted
a lot of interest among group theorists:

\begin{quote}
given $\varphi \in End(G)$, $w, t \in G$, find out whether or not
there is an element $s \in G$ such that $t = s^{-1} w \varphi(s)$.
\end{quote}

This problem is very non-trivial even for free groups; see
\cite{V} for an astonishing solution in the special
case where $\varphi$ is an {\it automorphism} of a free group. To
the best of our knowledge, this  decision problem is open for free
groups if $\varphi$ is an arbitrary endomorphism. Another class of
groups where the twisted conjugacy problem was considered is the
class of polycyclic-by-finite groups \cite{Troitsky}. Again, the
problem was solved for these groups in the special case where
$\varphi$ is an automorphism.

The conjugacy problem is a special case of the  twisted conjugacy
problem, where $\varphi$ is the identity map. Now a natural question
is: what is the advantage of the more general (double) twisted
conjugacy search problem over the conjugacy search problem in the
context of an authentication scheme? The answer is: if the platform
(semi)group $G$ has ``a lot" of endomorphisms, then Alice (the
prover), who selects $\varphi, \psi$, $w$, and  $s$, has an
opportunity to select them in such a way that there are a lot of
cancelations between $\psi(s), w$, and $\varphi(s)$, thus rendering
length based attacks ineffective.

In this paper, we use the semigroup of all $2 \times 2$ matrices
over truncated one-variable polynomials over ${\mathbf F}_2$, the
field of two elements, as the platform. It may seem that the
platform necessarily has  to be a group since one should at least
have the element $s$ (see above) invertible. However, as we will see
in the next section, we do not really need the invertibility to make
our authentication protocol work; what we need is just {\it some}
antihomomorphism of $G$ into itself, i.e., a map  $\ast: G \to G$
such that $(ab)^\ast =  b^\ast a^\ast$ for any $a, b \in G$. Every
group has such an antihomomorphism; it takes every element to its
inverse. Every {\it semigroup of square matrices} has such an
antihomomorphism, too; it takes every matrix to its transpose. Some
(semi)groups have other special antihomomorphisms; for example, any
free (semi)group has an antihomomorphism that rewrites every element
``backwards", i.e., right-to-left.  Here we prefer to focus on
semigroups of matrices (over commutative rings) since we believe
that these have several features making them fit to be platforms of
various cryptographic protocols, see \cite{Shpilrain_hash} for a
more detailed discussion.

\section{The protocol}
\label{protocol}

In this section, we give a description of a single round of our
authentication protocol. As with the original Fiat-Shamir scheme,
this protocol has to be repeated  $k$ times if one wants to reduce
the probability of successful forgery to $\frac{1}{2^k}$.

Here Alice is the prover and Bob the verifier. Let $G$ be the
platform semigroup, and  $\ast$  an antihomomorphism of $G$, i.e.,
$(ab)^\ast = b^\ast a^\ast$.

\begin{enumerate}
\item Alice's public key is a pair of endomorphisms $\varphi$, $\psi$ of the group $G$
and two elements $w, t \in G$, such that $t = \psi(s^\ast) w
\varphi(s)$, where $s \in G$ is her private key.

    \item To begin authentication,
Alice selects an element  $r \in G$ and sends the element $u =
\psi(r^\ast) t \varphi(r)$, called the {\em commitment}, to Bob.
    \item
Bob chooses a random bit $c$ and sends it to Alice.
\begin{itemize}
    \item
If $c=0$, then Alice sends $v = r$ to Bob and Bob checks if the
equality $u = \psi(v^\ast) t \varphi(v)$ is satisfied. If it is,
then Bob accepts the authentication.
    \item
If $c=1$, then Alice sends $v = sr$ to Bob and Bob checks if the
equality $u = \psi(v^\ast) w \varphi(v)$ is satisfied. If it is,
then Bob accepts the authentication.
\end{itemize}
\end{enumerate}

\noindent Let us check now that everything works the way we want it
to work.

\begin{itemize}
    \item
If $c=0$, then $v = r$, so $\psi(v^\ast) t \varphi(v) = \psi(r^\ast)
t \varphi(r) = u$.

\item
If $c=1$, then $v = sr$, so  $\psi(v^\ast) w \varphi(v) =
\psi((sr)^\ast) w \varphi(sr) = \psi(r^\ast s^\ast) w \varphi(s)
\varphi(r) = \psi(r^\ast) \psi(s^\ast) w \varphi(s) \varphi(r) = u$.

\end{itemize}

\section{The platform and parameters}
\label{platform}

Our suggested platform semigroup $G$ is the semigroup of all $2
\times 2$ matrices over truncated one-variable polynomials over
${\mathbf F}_2$, the field of two elements. Truncated (more
precisely,  $N$-truncated) one-variable polynomials over ${\mathbf
F}_2$ are expressions of the form ${\displaystyle \sum_{0 \le i \le
N-1} a_{i}x^i}$, where $a_{i}$ are elements of ${\mathbf F}_2$, and
$x$ is a variable. In other words, $N$-truncated
  polynomials are elements of the factor algebra of the
algebra ${\mathbf F}_2[x]$ of one-variable polynomials over
${\mathbf F}_2$ by the ideal generated by $x^N$.

Our semigroup $G$ has a lot of endomorphisms induced by
endomorphisms of the algebra of truncated polynomials. In fact, any
map of the form $x \to p(x)$, where  $p(x)$ is a truncated
polynomial with zero constant term, can be extended to an
endomorphism $\phi_p$ of the algebra of truncated polynomials.
Indeed, it is sufficient to show that $\phi_p(x^N)=(p(x))^N$ belongs
to the ideal generated by $x^N$, which is obviously the case if
$p(x)$ has zero constant term. Then, since $\phi_p$ is both an
additive and a multiplicative homomorphism, it extends to an
endomorphism of the semigroup of all $2 \times 2$ matrices over
truncated one-variable polynomials in the natural way.

If we now let the  antihomomorphism $\ast$ from the description of
the protocol in our Section \ref{protocol} to be the matrix
transposition, we have everything set up for an authentication
scheme using the semigroup $G$ as the platform.

Now we have to specify parameters involved in our scheme. The
parameter $N$ determines the size of the key space. If $N$ is on the
order of 300, then there are $2^{300}$ polynomials of degree $<N$
over ${\mathbf F}_2$, so there are $2^{1200}$ $2 \times 2$ matrices
over $N$-truncated polynomials, i.e., the size of the private key
space is $2^{1200}$, which is large enough.

At the same time, computations with (truncated) polynomials over
${\mathbf F}_2$ are very efficient (see e.g.  \cite{BP},
\cite{complexity}, or \cite{GG}  for details). In particular,
\begin{itemize}
    \item
Addition of two polynomials of  degree $N$ can be performed in
$O(N)$ time.
    \item
Multiplication of two polynomials of degree $N$ can be performed in $O(N \log_2N)$ time.
    \item
Computing composition $p(q(x)) \mod x^N$ of two polynomial of degree
$N$ can be performed in $O((N \log_2N)^{\frac{3}{2}})$ time (see
e.g. \cite[p.51]{complexity}).
\end{itemize}

Since those are the only operations used in our protocol, the time
complexity of executing a  single round of the protocol is $O((N
\log_2 N)^{\frac{3}{2}})$.

The size of public key space is large, too. One public key is,
again, a $2 \times 2$ matrix over $N$-truncated polynomials, and two
other public keys are endomorphisms of the form $x \to p(x)$, where
$p(x)$ is an $N$-truncated polynomial with zero constant term. Thus,
the number of different endomorphisms in this context is on the
order of $2^{300}$, hence the number of different {\it pairs} of
endomorphisms is on the order of $2^{600}$.

We also have to say a few words about how  a private key $s \in G$
is selected. We suggest that all  entries of the matrix $s$ have
non-zero constant term; other  coefficients of the entries can be
selected randomly, i.e., ``0"  and  ``1" are selected with
probability $\frac{1}{2}$ each. Non-zero constant terms are useful
here to ensure that there are sufficiently many non-zero terms in
the final product $t = \psi(s^\ast) w \varphi(s)$.

\section{Cryptanalysis}
\label{cryptanalysis}

As we have pointed out in  the previous section, the key space with
suggested parameters is quite large, so that a ``brute force" attack
by exhausting the key space is not feasible.

The next natural attack that comes to mind is attempting to solve a
system of equations over ${\mathbf F}_2$ that arises from equating
coefficients at the same powers of $x$ on both sides of the equation
$t = \psi(s^\ast) w \varphi(s)$.  Recall that in this equation $t,
w, \varphi,$ and $\psi$ are known, whereas $s$ is unknown.

More specifically, our experiments emulating this attack were
designed as follows. The entries of the private matrix $s$ were
generated as polynomials of degree $N-1$, with $N=100$ (which is
much smaller than the suggested $N=300$), with randomly selected
binary coefficients, except that the constant term in all
polynomials was 1. Then, the endomorphisms $\varphi$ and $\psi$ were
of the form $x \to p_i(x)$, where  $p_i(x)$ are polynomials of
degree $N-1$, with $N=150$, with randomly selected binary
coefficients, except that the constant term in both of them was 0.
Finally, the  entries of the public matrix $w$ were generated,
again,  as polynomials of degree $N-1$, with $N=100$, with randomly
selected binary coefficients, except that the constant term in all
polynomials was 1.

The attack itself then proceeds as follows. The matrix equation $t =
\psi(s^\ast) w \varphi(s)$ is converted to a system of $4N$
polynomial equations ($N$ for   each entry of a $2 \times 2$ matrix)
over ${\mathbf F}_2$. The unknowns in this system are coefficients
of the polynomials of degree $N-1$ that are the entries of the
private matrix $s$. Then, starting with the constant term and going
up, we equate coefficients at the same powers of $x$ on both sides
of each equation. After that, again starting with the coefficients
at the constant term and going up, we find all possible solutions of
each equation, one at a time. Thus we are getting a ``tree" of
solutions because some of the unknowns that occur in coefficients at
lower powers of $x$ also occur in coefficients at higher powers of
$x$. If this tree does not grow too fast, then there is a chance
that we can get all the way to the coefficients at highest power of
$x$, thereby finding a solution of the system. This solution may not
necessarily yield the same matrix $s$ that was selected by Alice,
but it is sufficient for forgery anyway.

We have run over 1000 experiments of this kind (which took about two
weeks),  allowing the solution tree to grow up to the width of
16384, i.e., allowing to go over at most 16384 solutions of each
equation when proceeding to a higher power of $x$. Each experiment
ran on a personal computer with Pentium 2Ghz dual core processor.
The success rate of the described attack with these parameters was
$0\%$.

\section{Conclusions}
\label{conclusions}

We have introduced:

\begin{enumerate}

\item An authentication scheme based on the (double) twisted conjugacy
problem, a new problem, which is allegedly hard in some
(semi)groups.

\item A new platform semigroup, namely the  semigroup of all $2
\times 2$ matrices over truncated one-variable polynomials over
${\mathbf F}_2$. Computation in this semigroup is very efficient
and, at the same time, the non-commutative structure of this
semigroup provides for security at least against obvious attacks.

\end{enumerate}

We point out here one important advantage of using the (double)
twisted conjugacy problem over using a more ``traditional" conjugacy
search problem as far as (semi)groups of matrices are concerned. The
conjugacy search problem admits a linear algebra attack upon
rewriting the equation $x^{-1}gx=h$ as  $gx=xh$; the latter
translates into a system of $n^2$ linear equations with $n^2$
unknowns, where $n$ is the size of the matrices involved, and the
unknowns are the entries of the matrix $x$.  Of course, if the
entries come not from a field but from a more general ring, such a
system of linear equations does not necessarily admit a
straightforward solution, but methods emulating standard techniques
(like Gauss elimination) usually have a pretty good success rate
anyway. For the twisted conjugacy problem, however, there is no
reduction to a system of linear equations.

We have considered an attack based on reducing the twisted conjugacy
problem to a system of {\it polynomial} equations over ${\mathbf
F}_2$, but this attack becomes computationally infeasible even with
a much smaller crucial parameter (which is the maximum degree of the
polynomials involved) than the one we suggest in this paper.

\end{document}